\documentclass{article}
\usepackage{graphicx} 
\usepackage[colorlinks]{hyperref}
\usepackage[all]{xy}
\usepackage{pdflscape}
\usepackage[hypcap=false]{caption}
\usepackage{amsmath,amssymb,amsfonts,latexsym,float,graphics,epsfig,amsthm} 
\usepackage{setspace}
\usepackage{xcolor}
\usepackage[top=2cm, bottom=2cm, left=2cm, right=2cm]{geometry} 

\newtheorem{theorem}{Theorem}[section]
\newtheorem{proposition}[theorem]{Proposition}
\newtheorem{lemma}[theorem]{Lemma}

\newtheorem{example}[theorem]{Example}

\title{Couplings of 3-anchored Bundles} 

\begin{document}

\maketitle
\begin{center} 
Begüm Ateşli$^{\ast,\ast\ast,}$\footnote{e-mails: 
\href{mailto:b.atesli@gtu.edu.tr}{b.atesli@gtu.edu.tr}, \href{mailto:begumatesli@itu.edu.tr}{begumatesli@itu.edu.tr} corresponding author,}, Oğul Esen $^{\ast,\dagger,}$\footnote{e-mail: 
\href{mailto:oesen@gtu.edu.tr}{oesen@gtu.edu.tr},} and Serkan Sütlü$^{\ast,}$\footnote{e-mail: 
\href{mailto:serkansutlu@gtu.edu.tr}{serkansutlu@gtu.edu.tr},}\\

\bigskip
$^\ast$Department of Mathematics, \\ Gebze Technical University, 41400 Gebze,
Kocaeli, Turkey.
\bigskip

$^{\ast\ast}$Department of Mathematics Engineering \\ İstanbul Technical University,  34467 Maslak, İstanbul

\bigskip

$^\dagger$Center for Mathematics and its Applications,\\  Khazar University, Baku, AZ1096,
Azerbaijan

\begin{abstract}

This work develops an algebraic framework for merging two $3$-anchored bundles
over the same base manifold, equipped with mutual actions and two twisted
cocycle terms, so as to obtain a $3$-Lie algebroid structure on the
corresponding Whitney sum. We also record the purely algebraic counterpart
of this construction, namely the bicocycle double cross product $3$-Lie
algebra, obtained by removing the anchor and Leibniz-type compatibility
conditions. The resulting framework provides a unified setting for $3$-Lie
algebroids and contains, as special cases, unified products, double cross
products, semi-direct products, cocycle extensions, and direct products.

\smallskip

\smallskip


\end{abstract}

\end{center}

\tableofcontents
\onehalfspacing

\setlength{\parindent}{2em}
\setlength{\parskip}{3ex}

 \section{Introduction}

Lie algebroids constitute a natural and powerful geometric framework that simultaneously extends the notions of Lie algebras and tangent bundles \cite{MackenzieDG,Mackenzie-book,Para67}. Owing to this unifying character, they provide a flexible language for both algebraic and geometric analysis and have found substantial applications in geometric mechanics and dynamical systems; see, for instance, \cite{GrabowskaAlg,marle2014lie,martinez2001,WeinLag}.

A higher-order extension of this framework is given by $n$-Lie algebroids, also referred to as Filippov algebroids \cite{GrabowskiMarmo2000,Mishra,Vallejo}. These structures arise by replacing the binary Lie bracket with an $n$-ary Lie bracket. Accordingly, the classical Jacobi identity, which encodes the compatibility and integrability properties of the binary bracket, is replaced by the fundamental identity, also known as the Takhtajan--Filippov identity. In the present work, we focus on the case $n=3$ and study coupling constructions for $3$-Lie algebroids from the viewpoint of anchored bundles.

\noindent \textbf{(De)Coupling Problem.}
The problem addressed in the present paper is the ``(de)coupling problem''. This problem concerns the construction of a larger algebraic or geometric structure from two given constituent structures in such a way that the original components remain visible in the coupled object. The theory of (de)coupling has been studied in the literature in various contexts. The coupling of a Lie group or a Lie algebra with its representation space leads to the semidirect product theory \cite{CeHoMaRa98, kuper83, marsden83b, marsden1984semidirect}, which has various applications in rigid body dynamics, as well as in fluid and plasma theories. A generalization of this construction is the double cross product, or matched pair, theory, which allows mutual interactions between Lie groups or Lie algebras \cite{Ma90, majid1990physicsa}. This has been applied to Lagrangian and Hamiltonian dynamics \cite{esen2016hamiltonian, esen2017lagrangian}, higher-order dynamics \cite{EsenKudeSutlu21}, dissipative systems \cite{esen2021extensions}, and discrete dynamics \cite{esen2018matched}, with extensions to electromagnetic theory \cite{Esvd17}, plasma-fluid relations \cite{EsGrMiGu19, EsSu21}, and chemical kinetics \cite{AjChEsGrKlPa}. On the Lie algebroid level, this framework is studied in \cite{Mokr97}.  

A more general framework is provided by unified product theory, also called cocycle double cross product theory, in which twisted cocycles are incorporated to couple Lie algebras or Lie algebroids with vector bundles \cite{AgorMili11, Agore-Lie, AgorMili14-II, AgoreMilitaru-book, esen2021bicocycle}. This method allows one to describe and decompose algebroid structures in terms of subalgebroids and complementary vector bundles, and it has also been used in the study of dynamics on homogeneous spaces \cite{uccgun2024dynamics}. Recently, bicocycle double cross products (BDCPs) were introduced by allowing two twisted cocycles in the coupling of vector bundles into a Lie algebroid \cite{atecsli2025product,esen2021bicocycle}. The aim of the present paper is to develop an analogue of this BDCP strategy for $3$-Lie algebroids.

\noindent\textbf{Goal.}
The main goal of this work is to couple two $3$-anchored bundles over the same base manifold in order to construct a $3$-Lie algebroid structure on their Whitney sum. This is achieved through a bicocycle double cross product strategy, which allows mutual \textit{actions} together with two twisted cocycle terms. The resulting construction contains, as particular cases, semidirect product $3$-Lie algebroids, double cross product or matched pair $3$-Lie algebroids, unified product $3$-Lie algebroids, and cocycle extension $3$-Lie algebroids. In addition, we record the purely algebraic counterpart of the construction,
namely the bicocycle double cross product $3$-Lie algebra, which is obtained by
removing the anchor and Leibniz-type compatibility conditions from the
$3$-Lie algebroid setting.

The coupling problem also has a natural decoupling counterpart. Namely, once a $3$-Lie algebroid is given together with a decomposition of its underlying vector bundle, one may ask how the original $3$-Lie algebroid structure is encoded by brackets, actions, and cocycle terms on the components. We shall therefore present both the coupling and decoupling aspects from a purely algebraic standpoint, deriving the compatibility conditions for the former and formulating the corresponding structural interpretation for the latter.

\section{3-Lie Algebroids} \label{coupling}

Consider a vector bundle $\tau:\mathcal{A}\mapsto M$. A $3$-anchor is defined to be a map from the skew product $\Gamma(\mathcal{A})\wedge\Gamma(\mathcal{A})$ of the sections of $\mathcal{A}$ to the space $\Gamma(TM)$ of vector fields on $M$:
\begin{equation}
\mathfrak{a}_{\mathcal{A}}:\Gamma(\mathcal{A})\wedge\Gamma(\mathcal{A}) \longrightarrow \Gamma(TM).
\end{equation}
In this case, the quadruple $(\mathcal{A},\tau,M,\mathfrak{a}_\mathcal{A})$ is said to be a $3$-anchored bundle. 

\noindent\textbf{$3$-Lie Algebroids.}
A $3$-Lie algebroid, also called a $3$-Filippov algebroid, see \cite{GrabowskiMarmo2000,Mishra,Vallejo}, is a $3$-anchored bundle $(\mathcal{A},\tau,M,\mathfrak{a}_\mathcal{A})$ equipped with a $3$-Lie bracket operation 
\begin{equation}\label{eq-Filippov-3-bracket}
[\bullet,\bullet,\bullet]:\wedge^3 \Gamma(\mathcal{A}) \longrightarrow \Gamma(\mathcal{A})
\end{equation}
so that the $3$-anchor is a Leibniz algebra map given by
\begin{equation}\label{3-Lie-cond1}
\begin{split} 
 \mathfrak{a}_{\mathcal{A}}([X_1\wedge X_2, Y_1\wedge Y_2]) 
 =
 [\mathfrak{a}_{\mathcal{A}}(X_1\wedge X_2),
 \mathfrak{a}_{\mathcal{A}}(Y_1\wedge Y_2)] .
\end{split}
\end{equation}
Here, the bracket on the right hand side is the Jacobi--Lie bracket of vector fields on $M$, whereas the bracket on the left hand side is the Leibniz algebra bracket on $\wedge^2\Gamma(\mathcal{A})$ defined by
\begin{equation}\label{SS-bracket-3}
[X_1\wedge X_2, Y_1\wedge Y_2]
:=
[X_1,X_2,Y_1]\wedge Y_2
+
Y_1\wedge [X_1,X_2,Y_2],
\end{equation}
for any $X_1,X_2,Y_1,Y_2$ in $\Gamma(\mathcal{A})$. 

To have a $3$-Lie algebroid structure, we further ask that the fundamental identity, also called the Filippov or Takhtajan identity,
\begin{equation}\label{3-Lie-cond2}
[X_1,X_2,[Y_1,Y_2,Y_3]]
=
[[X_1,X_2,Y_1],Y_2,Y_3]
+
[Y_1,[X_1,X_2,Y_2],Y_3]
+
[Y_1,Y_2,[X_1,X_2,Y_3]],
\end{equation}
and the Leibniz rule
\begin{equation}\label{3-Lie-cond3} 
[X_1,X_2,fY]
=
f[X_1,X_2,Y]
+
\mathfrak{a}_{\mathcal{A}}(X_1\wedge X_2)(f)Y
\end{equation}
are both satisfied for any $X_1,X_2,Y_1,Y_2,Y_3,Y$ in $\Gamma(\mathcal{A})$, and any $f$ in $C^\infty(M)$. Accordingly, we denote a $3$-Lie algebroid by a quintuple 
\begin{equation}
(\mathcal{A},\tau,M,\mathfrak{a}_\mathcal{A},[\bullet,\bullet,\bullet]). 
\end{equation}

\noindent\textbf{$3$-Lie Algebroid Actions on Vector Bundles.}
Let 
$(\mathcal{A},\tau,M,\mathfrak{a}_{\mathcal{A}},[\bullet,\bullet,\bullet])$ 
be a $3$-Lie algebroid, and let $\pi:E\mapsto M$ be a vector bundle.
A $3$-Lie algebroid action, or representation, of $\mathcal{A}$ on $E$ is an
$\mathbb{R}$-bilinear skew-symmetric map
\begin{equation}\label{3-Lie-algebroid-action-on-vector-bundle}
\varrho:
\Gamma(\mathcal{A})\wedge \Gamma(\mathcal{A})
\longrightarrow
\operatorname{End}_{\mathbb{R}}(\Gamma(E)),
\qquad
X_1\wedge X_2
\mapsto
\varrho(X_1,X_2).
\end{equation}
Here, $\operatorname{End}_{\mathbb{R}}(\Gamma(E))$ denotes the space of
$\mathbb{R}$-linear endomorphisms of $\Gamma(E)$. In particular, the operators
$\varrho(X_1,X_2)$ are not assumed to be $C^\infty(M)$-linear. Their failure
to be $C^\infty(M)$-linear is controlled by the anchor through the Leibniz rule
below.

The defining conditions naturally split into two groups. The first group
expresses the compatibility of the action with the $C^\infty(M)$-module
structures. More precisely, the action is $C^\infty(M)$-linear in the
$\mathcal{A}$-variables, while its action on sections of $E$ is governed by
the anchor map:
\begin{equation}\label{3-Lie-action-Cinfty}
\begin{split} 
\varrho(fX_1,X_2)(s)
&=
f\varrho(X_1,X_2)(s),\\
\varrho(X_1,X_2)(fs)
&=
f\varrho(X_1,X_2)(s)
+
\mathfrak{a}_{\mathcal{A}}(X_1\wedge X_2)(f)s.
\end{split}
\end{equation}

The second group consists of the algebraic representation identities. These
conditions say that the operators $\varrho(X_1,X_2)$ represent the Leibniz
algebra structure induced on $\wedge^2\Gamma(\mathcal{A})$ by the $3$-Lie
bracket, and that this representation is compatible with the fundamental
identity of the $3$-Lie algebroid:
\begin{equation}\label{3-Lie-action-representation}
\begin{split}
&
\varrho(X_1,X_2)\varrho(Y_1,Y_2)
-
\varrho(Y_1,Y_2)\varrho(X_1,X_2)
\\
&\qquad =
\varrho([X_1,X_2,Y_1],Y_2)
+
\varrho(Y_1,[X_1,X_2,Y_2]),
\\
&
\varrho([Y_1,Y_2,Y_3],Y_4)
\\
&\qquad =
\varrho(Y_1,Y_2)\varrho(Y_3,Y_4)
+
\varrho(Y_2,Y_3)\varrho(Y_1,Y_4)
+
\varrho(Y_3,Y_1)\varrho(Y_2,Y_4),
\end{split}
\end{equation}
for any $X_1,X_2,Y_1,Y_2,Y_3,Y_4$ in $\Gamma(\mathcal{A})$,
$s$ in $\Gamma(E)$, and $f$ in $C^\infty(M)$. In this case, the vector bundle $E$ is called a representation bundle of the
$3$-Lie algebroid $\mathcal{A}$.

\section{Couplings of 3-anchored Bundles}
     
We assume two $3$-anchored bundles
$(\mathcal{A},\tau,M,\mathfrak{a}_{\mathcal{A}})$ and
$(\mathcal{B},\kappa,M,\mathfrak{a}_{\mathcal{B}})$ over the same base
manifold $M$, with $3$-anchor maps
\begin{equation}
\mathfrak{a}_{\mathcal{A}}:
\Gamma(\mathcal{A})\wedge\Gamma(\mathcal{A})
\longrightarrow
\Gamma(TM),
\qquad
\mathfrak{a}_{\mathcal{B}}:
\Gamma(\mathcal{B})\wedge\Gamma(\mathcal{B})
\longrightarrow
\Gamma(TM),
\end{equation}
respectively. We denote by $\mathcal{A}\times\mathcal{B}$ the Whitney sum, that is, the fiberwise direct product, of $\mathcal{A}$ and $\mathcal{B}$ over $M$. Thus,
an element of the fiber over $x\in M$ is a pair
\begin{equation}
(X_x,\eta_x)\in \mathcal{A}_x\times\mathcal{B}_x,
\end{equation}
where $X_x\in \mathcal{A}_x$ and $\eta_x\in \mathcal{B}_x$. Accordingly, a
section of $\mathcal{A}\times\mathcal{B}$ is written as a pair
\begin{equation}
(X,\eta)\in \Gamma(\mathcal{A}\times\mathcal{B}),
\end{equation}
where $X\in \Gamma(\mathcal{A})$ and $\eta\in \Gamma(\mathcal{B})$. Define
\begin{equation}
\begin{split}
&\mathfrak{a}_{\mathcal{A}\times\mathcal{B}}:
\Gamma(\mathcal{A}\times\mathcal{B})
\wedge
\Gamma(\mathcal{A}\times\mathcal{B})
\longrightarrow
\Gamma(TM), \\
&\mathfrak{a}_{\mathcal{A}\times\mathcal{B}}
\big((X_1,\eta_1)\wedge(X_2,\eta_2)\big)
=
\mathfrak{a}_{\mathcal{A}}(X_1\wedge X_2)
+
\mathfrak{a}_{\mathcal{B}}(\eta_1\wedge\eta_2).
\end{split}
\end{equation}
Then
\begin{equation}
(\mathcal{A}\times\mathcal{B},
\tau\times\kappa,
M,
\mathfrak{a}_{\mathcal{A}\times\mathcal{B}})
\end{equation}
is a $3$-anchored bundle.

\noindent\textbf{$3$-Bracket on the Whitney sum.}
We are now ready to study the algebraic structure on the Whitney sum
$\mathcal{A}\times\mathcal{B}$. For this purpose, we first assume the following
$\mathbb{R}$-trilinear skew-symmetric operations:
\begin{equation}\label{sennur2}
\phi:\wedge^3\Gamma(\mathcal{A})\longrightarrow\Gamma(\mathcal{A}),
\qquad
\zeta:\wedge^3\Gamma(\mathcal{A})\longrightarrow\Gamma(\mathcal{B}).
\end{equation}
These two maps define an $\mathcal{A}\times\mathcal{B}$-valued $3$-\textit{bracket}
on $\Gamma(\mathcal{A})$ by
\begin{equation}\label{3-Lie-A}
[\bullet,\bullet,\bullet]_{\mathcal{A}}:
\wedge^3\Gamma(\mathcal{A})
\longrightarrow
\Gamma(\mathcal{A}\times\mathcal{B}),
\qquad
[X_1,X_2,X_3]_{\mathcal{A}}
=
\big(
\phi(X_1,X_2,X_3),
\zeta(X_1,X_2,X_3)
\big),
\end{equation}
for all $X_1,X_2,X_3$ in $\Gamma(\mathcal{A})$. The existence of a non-zero
$\zeta$ implies that $\mathcal{A}$ is not closed under the operation
$[\bullet,\bullet,\bullet]_{\mathcal{A}}$. Consequently,
$[\bullet,\bullet,\bullet]_{\mathcal{A}}$ cannot be regarded as a genuine
$3$-bracket on $\mathcal{A}$. It becomes a closed operation, and hence a genuine
$3$-bracket on $\mathcal{A}$, if and only if $\zeta$ is identically zero. For
this reason, the map $\zeta$ will be called a twisted cocycle.

Similarly, for the $3$-anchored bundle $\mathcal{B}$, we assume the following
$\mathbb{R}$-trilinear skew-symmetric operations:
\begin{equation}\label{sennur1}
\psi:\wedge^3\Gamma(\mathcal{B})\longrightarrow\Gamma(\mathcal{A}),
\qquad
\theta:\wedge^3\Gamma(\mathcal{B})\longrightarrow\Gamma(\mathcal{B}).
\end{equation}
These two maps define an $\mathcal{A}\times\mathcal{B}$-valued $3$-\textit{bracket}
on $\Gamma(\mathcal{B})$ by
\begin{equation}\label{3-Lie-B}
[\bullet,\bullet,\bullet]_{\mathcal{B}}:
\wedge^3\Gamma(\mathcal{B})
\longrightarrow
\Gamma(\mathcal{A}\times\mathcal{B}),
\qquad
[\eta_1,\eta_2,\eta_3]_{\mathcal{B}}
=
\big(
\psi(\eta_1,\eta_2,\eta_3),
\theta(\eta_1,\eta_2,\eta_3)
\big),
\end{equation}
for all $\eta_1,\eta_2,\eta_3$ in $\Gamma(\mathcal{B})$. The existence of a
non-zero $\psi$ implies that $\mathcal{B}$ is not closed under the operation
$[\bullet,\bullet,\bullet]_{\mathcal{B}}$. Thus,
$[\bullet,\bullet,\bullet]_{\mathcal{B}}$ cannot be regarded as a genuine
$3$-bracket on $\mathcal{B}$. It becomes a closed operation, and hence a genuine
$3$-bracket on $\mathcal{B}$, if and only if $\psi$ is identically zero. For
this reason, the map $\psi$ will also be called a twisted cocycle.

We further consider the following $\mathbb{R}$-trilinear maps:
\begin{equation}\label{sennur3}
\rho:
\Gamma(\mathcal{B})\wedge\Gamma(\mathcal{B})
\otimes\Gamma(\mathcal{A})
\longrightarrow
\Gamma(\mathcal{A}),
\qquad
\sigma:
\Gamma(\mathcal{B})
\otimes
\Gamma(\mathcal{A})\wedge\Gamma(\mathcal{A})
\longrightarrow
\Gamma(\mathcal{B}).
\end{equation}
At first glance, $\rho$ may be viewed as a left \textit{action} of
$\mathcal{B}$ on $\mathcal{A}$, while $\sigma$ may be viewed as a right
\textit{action} of $\mathcal{A}$ on $\mathcal{B}$. However, at this level they
are only action-type terms; in general, they need not satisfy the identities
required of genuine $3$-Lie algebroid actions. In special cases, in particular
when the twisted cocycles $\zeta$ and $\psi$ vanish, these maps become genuine
$3$-Lie algebroid actions.

Using the algebraic data
\begin{equation}
(\phi,\zeta,\psi,\theta,\rho,\sigma),
\end{equation}
we define a $3$-bracket on $\Gamma(\mathcal{A}\times\mathcal{B})$ by
\begin{equation}\label{BDCP-3-Lie-bra}
\begin{split}
&[\bullet,\bullet,\bullet]_{{}_{\zeta}\bowtie_{\psi}}:
\wedge^3\Gamma(\mathcal{A}\times\mathcal{B})
\longrightarrow
\Gamma(\mathcal{A}\times\mathcal{B}),
\\
&[(X_1,\eta_1),(X_2,\eta_2),(X_3,\eta_3)]_{{}_{\zeta}\bowtie_{\psi}}
=
(\widetilde{X},\widetilde{\eta}),
\end{split}
\end{equation}
where
\begin{equation}\label{BDCP-3-Lie-bra+}
\begin{split}
\widetilde{X}
&=
\phi(X_1,X_2,X_3)
+
\psi(\eta_1,\eta_2,\eta_3)
+
\rho(\eta_2,\eta_3,X_1)
-
\rho(\eta_1,\eta_3,X_2)
+
\rho(\eta_1,\eta_2,X_3),
\\
\widetilde{\eta}
&=
\theta(\eta_1,\eta_2,\eta_3)
+
\zeta(X_1,X_2,X_3)
+
\sigma(\eta_3,X_1,X_2)
-
\sigma(\eta_2,X_1,X_3)
+
\sigma(\eta_1,X_2,X_3).
\end{split}
\end{equation}

\noindent\textbf{$3$-Lie Algebroid Structure on the Whitney sum.}
To establish that $\mathcal{A}\times\mathcal{B}$ is a $3$-Lie algebroid, the conditions in \eqref{3-Lie-cond1}, \eqref{3-Lie-cond2} and \eqref{3-Lie-cond3} must be satisfied. To this end, we derive below the corresponding list of compatibility conditions.

The first set of compatibility conditions arises from \eqref{3-Lie-cond1}, which states that the anchor map is a Leibniz algebra morphism: 
\begin{equation}\label{doga1}
    \begin{split}
        &
        [\mathfrak{a}_{\mathcal{A}}(X_1\wedge X_2),
        \mathfrak{a}_{\mathcal{A}}(X_3\wedge X_4)]
        =
        \mathfrak{a}_{\mathcal{A}}(\phi(X_1,X_2,X_3)\wedge X_4)
        +
        \mathfrak{a}_{\mathcal{A}}(X_3\wedge \phi(X_1,X_2,X_4)),
        \\
        &
        \mathfrak{a}_{\mathcal{B}}(\zeta(X_1,X_2,X_3)\wedge\eta_4)
        =
        0,
        \\
        &
        [\mathfrak{a}_{\mathcal{A}}(X_1\wedge X_2),
        \mathfrak{a}_{\mathcal{B}}(\eta_3\wedge\eta_4)]
        =
        \mathfrak{a}_{\mathcal{B}}(\sigma(\eta_3,X_1,X_2)\wedge\eta_4)
        +
        \mathfrak{a}_{\mathcal{B}}(\eta_3\wedge\sigma(\eta_4,X_1,X_2)),
        \\
        &
        \mathfrak{a}_{\mathcal{B}}(\sigma(\eta_2,X_1,X_3)\wedge\eta_4)
        -
        \mathfrak{a}_{\mathcal{A}}(X_3\wedge\rho(\eta_2,\eta_4,X_1))
        =
        0,
        \\
        &
        [\mathfrak{a}_{\mathcal{B}}(\eta_1\wedge\eta_2),
        \mathfrak{a}_{\mathcal{A}}(X_3\wedge X_4)]
        =
        \mathfrak{a}_{\mathcal{A}}(\rho(\eta_1,\eta_2,X_3)\wedge X_4)
        +
        \mathfrak{a}_{\mathcal{A}}(X_3\wedge\rho(\eta_1,\eta_2,X_4)),
        \\
        &
        \mathfrak{a}_{\mathcal{A}}(X_3\wedge\psi(\eta_1,\eta_2,\eta_4))
        =
        0,
        \\
        &
        [\mathfrak{a}_{\mathcal{B}}(\eta_1\wedge\eta_2),
        \mathfrak{a}_{\mathcal{B}}(\eta_3\wedge\eta_4)]
        =
        \mathfrak{a}_{\mathcal{B}}(\theta(\eta_1,\eta_2,\eta_3)\wedge\eta_4)
        +
        \mathfrak{a}_{\mathcal{B}}(\eta_3\wedge\theta(\eta_1,\eta_2,\eta_4)).
    \end{split}
\end{equation}

We next turn to the compatibility conditions arising from the fundamental identity \eqref{3-Lie-cond2}. The following conditions are obtained by applying the fundamental identity either to five sections of $\mathcal{A}$ or to five sections of $\mathcal{B}$, and then comparing the $\mathcal{A}$- and $\mathcal{B}$-components separately:
\begin{equation}\label{pelin1}
    \begin{split}
        &\phi(X_1,X_2,\phi(X_3,X_4,X_5)) = \phi(\phi(X_1,X_2,X_3),X_4,X_5) \\
        & \qquad \qquad + \phi(X_3,\phi(X_1,X_2,X_4),X_5) + \phi(X_3,X_4,\phi(X_1,X_2,X_5)), \\
        &\zeta(X_1,X_2,\phi(X_3,X_4,X_5)) + \sigma(\zeta(X_3,X_4,X_5),X_1,X_2) = \\ &\hspace{2cm}\zeta(\phi(X_1,X_2,X_3),X_4,X_5)  + \sigma(\zeta(X_1,X_2,X_3),X_4,X_5) + \\ 
        & \zeta(X_3,\phi(X_1,X_2,X_4),X_5) - \sigma(\zeta(X_1,X_2,X_4),X_3,X_5)+ \\
        & \qquad \qquad  \zeta(X_3,X_4,\phi(X_1,X_2,X_5)) + \sigma(\zeta(X_1,X_2,X_5),X_3,X_4), \\
        & \theta(\eta_1,\eta_2,\theta(\eta_3,\eta_4,\eta_5)) = \theta(\theta(\eta_1,\eta_2,\eta_3),\eta_4,\eta_5) + \\
        & \qquad \qquad  \theta(\eta_3,\theta(\eta_1,\eta_2,\eta_4),\eta_5) + \theta(\eta_3,\eta_4,\theta(\eta_1,\eta_2,\eta_5)), \\
         &\psi(\eta_1,\eta_2,\theta(\eta_3,\eta_4,\eta_5)) + \rho(\eta_1,\eta_2,\psi(\eta_3,\eta_4,\eta_5)) =  \\
        & \hspace{1cm}\psi(\theta(\eta_1,\eta_2,\eta_3),\eta_4,\eta_5) + \rho(\eta_4,\eta_5,\psi(\eta_1,\eta_2,\eta_3)) + \\ 
        &\hspace{2cm} \psi(\eta_3,\theta(\eta_1,\eta_2,\eta_4),\eta_5) - \rho(\eta_3,\eta_5,\psi(\eta_1,\eta_2,\eta_4)) + \\
        & \hspace{3cm}  \psi(\eta_3,\eta_4,\theta(\eta_1,\eta_2,\eta_5)) + \rho(\eta_3,\eta_4,\psi(\eta_1,\eta_2,\eta_5)).
    \end{split}
\end{equation}

The following conditions are obtained by applying the
fundamental identity to terms involving four sections of $\mathcal{A}$ and one
section of $\mathcal{B}$, and then comparing the $\mathcal{A}$- and
$\mathcal{B}$-components separately:
\begin{equation}\label{pelin2}
    \begin{split}
     &\sigma(\sigma(\eta_5,X_3,X_4),X_1,X_2) 
     =
     \sigma(\eta_5,\phi(X_1,X_2,X_3),X_4) \\
     & \qquad \qquad
     + \sigma(\eta_5,X_3,\phi(X_1,X_2,X_4)) 
     + \sigma(\sigma(\eta_5,X_1,X_2),X_3,X_4),
     \\
     &\rho(\zeta(X_1,X_2,X_3),\eta_5,X_4) 
     -
     \rho(\zeta(X_1,X_2,X_4),\eta_5,X_3)
     =
     0,
     \\
     &\sigma(\eta_2,X_1,\phi(X_3,X_4,X_5)) 
     =
     \sigma(\sigma(\eta_2,X_1,X_3),X_4,X_5) \\
     & \qquad \qquad
     - \sigma(\sigma(\eta_2,X_1,X_4),X_3,X_5) 
     + \sigma(\sigma(\eta_2,X_1,X_5),X_3,X_4),
     \\
     &\rho(\eta_2,\zeta(X_3,X_4,X_5),X_1)
     =
     0.
    \end{split}
\end{equation}

The following conditions are those obtained from the fundamental identity by
considering terms involving three sections of $\mathcal{A}$ and two sections of
$\mathcal{B}$, and then comparing the $\mathcal{A}$- and $\mathcal{B}$-components
separately:
\begin{equation}\label{pelin3}
    \begin{split} 
        &\phi(X_1,X_2,\rho(\eta_4,\eta_5,X_3)) 
        =
        \psi(\zeta(X_1,X_2,X_3),\eta_4,\eta_5) 
        + \rho(\eta_4,\eta_5,\phi(X_1,X_2,X_3)) \\
        & \qquad \qquad
        + \rho(\sigma(\eta_4,X_1,X_2),\eta_5,X_3) 
        + \rho(\eta_4,\sigma(\eta_5,X_1,X_2),X_3),
        \\
        &\zeta(X_1,X_2,\rho(\eta_4,\eta_5,X_3)) 
        =
        \theta(\zeta(X_1,X_2,X_3),\eta_4,\eta_5),
        \\
        &\rho(\eta_2,\sigma(\eta_5,X_3,X_4),X_1) 
        =
        \rho(\sigma(\eta_2,X_1,X_3),\eta_5,X_4) \\
        & \qquad \qquad
        - \rho(\sigma(\eta_2,X_1,X_4),\eta_5,X_3) 
        + \phi(X_3,X_4,\rho(\eta_2,\eta_5,X_1)),
        \\
        &\zeta(X_3,X_4,\rho(\eta_2,\eta_5,X_1)) 
        =
        0,
        \\
        &\rho(\eta_1,\eta_2,\phi(X_3,X_4,X_5)) 
        + \psi(\eta_1,\eta_2,\zeta(X_3,X_4,X_5)) \\
        &\qquad =
        \phi(\rho(\eta_1,\eta_2,X_3),X_4,X_5)
        + \phi(X_3,\rho(\eta_1,\eta_2,X_4),X_5) \\
        & \qquad \qquad
        + \phi(X_3,X_4,\rho(\eta_1,\eta_2,X_5)),
        \\
        &\theta(\eta_1,\eta_2,\zeta(X_3,X_4,X_5)) 
        =
        \zeta(\rho(\eta_1,\eta_2,X_3),X_4,X_5) \\
        & \qquad \qquad
        + \zeta(X_3,\rho(\eta_1,\eta_2,X_4),X_5) 
        + \zeta(X_3,X_4,\rho(\eta_1,\eta_2,X_5)).
    \end{split}
\end{equation}

The following conditions are those obtained from the fundamental identity by
considering terms involving two sections of $\mathcal{A}$ and three sections of
$\mathcal{B}$, and then comparing the $\mathcal{A}$- and $\mathcal{B}$-components
separately:
\begin{equation}\label{pelin4}
    \begin{split}
        &\phi(X_1,X_2,\psi(\eta_3,\eta_4,\eta_5)) = \psi(\sigma(\eta_3,X_1,X_2),\eta_4,\eta_5) \\
        & \qquad \qquad + \psi(\eta_3,\sigma(\eta_4,X_1,X_2),\eta_5) + \psi(\eta_3,\eta_4,\sigma(\eta_5,X_1,X_2)), \\
        & \zeta(X_1,X_2,\psi(\eta_3,\eta_4,\eta_5)) + \sigma(\theta(\eta_3,\eta_4,\eta_5),X_1,X_2) = \theta(\sigma(\eta_3,X_1,X_2),\eta_4,\eta_5) \\
        & \qquad \qquad + \theta(\eta_3,\sigma(\eta_4,X_1,X_2),\eta_5) + \theta(\eta_3,\eta_4,\sigma(\eta_5,X_1,X_2)) \\
        & \sigma(\eta_2,X_1,\rho(\eta_4,\eta_5,X_3)) = \theta(\sigma(\eta_2,X_1,X_3),\eta_4,\eta_5) \\
        & \qquad \qquad -\sigma(\eta_5,X_3,\rho(\eta_2,\eta_4,X_1)) + \sigma(\eta_4,X_3,\rho(\eta_2,\eta_5,X_1)), \\
        & \psi(\sigma(\eta_2,X_1,X_3),\eta_4,\eta_5) = 0, \\
        & \psi(\eta_1,\eta_2,\sigma(\eta_5,X_3,X_4)) = \phi(X_3,X_4,\psi(\eta_1,\eta_2,\eta_5)), \\
        & \theta(\eta_1,\eta_2,\sigma(\eta_5,X_3,X_4)) = \sigma(\eta_5,\rho(\eta_1,\eta_2,X_3),X_4) + \sigma(\eta_5,X_3,\rho(\eta_1,\eta_2,X_4)) \\
        & \qquad \qquad + \zeta(X_3,X_4,\psi(\eta_1,\eta_2,\eta_5)) + \sigma(\theta(\eta_1,\eta_2,\eta_5),X_3,X_4).
        \end{split}
        \end{equation}

The following conditions are those obtained from the fundamental identity by
considering terms involving one section of $\mathcal{A}$ and four sections of
$\mathcal{B}$, and then comparing the $\mathcal{A}$- and $\mathcal{B}$-components
separately:
\begin{equation}\label{pelin5}
\begin{split}
&\rho(\eta_2,\theta(\eta_3,\eta_4,\eta_5),X_1)
=
\rho(\eta_4,\eta_5,\rho(\eta_2,\eta_3,X_1)) \\
&\qquad \qquad
-
\rho(\eta_3,\eta_5,\rho(\eta_2,\eta_4,X_1))
+
\rho(\eta_3,\eta_4,\rho(\eta_2,\eta_5,X_1)),
\\
&\sigma(\eta_2,\psi(\eta_3,\eta_4,\eta_5),X_1)
=
0,
\\
&\rho(\eta_1,\eta_2,\rho(\eta_4,\eta_5,X_3))
=
\rho(\eta_4,\eta_5,\rho(\eta_1,\eta_2,X_3)) \\
&\qquad \qquad
+
\rho(\theta(\eta_1,\eta_2,\eta_4),\eta_5,X_3)
+
\rho(\eta_4,\theta(\eta_1,\eta_2,\eta_5),X_3),
\\
&\sigma(\eta_5,X_3,\psi(\eta_1,\eta_2,\eta_4))
-
\sigma(\eta_4,X_3,\psi(\eta_1,\eta_2,\eta_5))
=
0.
\end{split}
\end{equation}

We conclude the list of compatibility conditions with those arising from the
Leibniz rule \eqref{3-Lie-cond3}. Accordingly, one obtains:
\begin{equation}\label{doga2}
    \begin{split}
        & \phi(X_1,X_2,fX) 
        =
        f\phi(X_1,X_2,X) 
        + \mathfrak{a}_{\mathcal{A}}(X_1\wedge X_2)(f)X,
        \\
        & \zeta(X_1,X_2,fX) 
        =
        f\zeta(X_1,X_2,X),
        \\
        & \sigma(f\eta,X_1,X_2) 
        =
        f\sigma(\eta,X_1,X_2) 
        + \mathfrak{a}_{\mathcal{A}}(X_1\wedge X_2)(f)\eta,
        \\
        & \sigma(\eta_2,X_1,fX) 
        =
        f\sigma(\eta_2,X_1,X),
        \\
        & \rho(f\eta,\eta_2,X_1) 
        =
        f\rho(\eta,\eta_2,X_1),
        \\
        & \rho(\eta_1,\eta_2,fX) 
        =
        f\rho(\eta_1,\eta_2,X) 
        + \mathfrak{a}_{\mathcal{B}}(\eta_1 \wedge \eta_2)(f)X,
        \\
        & \psi(\eta_1,\eta_2,f\eta) 
        =
        f\psi(\eta_1,\eta_2,\eta),
        \\
        & \theta(\eta_1,\eta_2,f\eta) 
        =
        f\theta(\eta_1,\eta_2,\eta) 
        + \mathfrak{a}_{\mathcal{B}}(\eta_1 \wedge \eta_2)(f)\eta.
    \end{split}
\end{equation}
Combining the anchor compatibility conditions \eqref{doga1}, the compatibility
conditions obtained from the fundamental identity
\eqref{pelin1}--\eqref{pelin5}, and the Leibniz-type conditions
\eqref{doga2}, we obtain the following coupling result. This result gives the
bicocycle double cross product construction of $3$-Lie algebroids.

\begin{proposition}\label{thm-muratcan-off+}
Let
$(\mathcal{A},\tau,M,\mathfrak{a}_{\mathcal{A}})$
and
$(\mathcal{B},\kappa,M,\mathfrak{a}_{\mathcal{B}})$
be two $3$-anchored bundles over the same base manifold $M$, equipped with the
maps introduced in \eqref{sennur2}, \eqref{sennur1}, and \eqref{sennur3}.
Then the Whitney sum $\mathcal{A}\times\mathcal{B}$ carries a
bicocycle double cross product $3$-Lie algebroid structure with anchor
$\mathfrak{a}_{\mathcal{A}\times\mathcal{B}}$ and bracket
\eqref{BDCP-3-Lie-bra}--\eqref{BDCP-3-Lie-bra+} if and only if the
compatibility conditions
\eqref{doga1}, \eqref{pelin1}--\eqref{pelin5}, and \eqref{doga2}
are all satisfied.
\end{proposition}

 \section{Unified Product, Double Cross and Semi-Direct Products}

Having constructed the bicocycle double cross product $3$-Lie algebroid in its
most general form, we now examine several distinguished special cases obtained
by imposing natural restrictions on the twisted cocycles and action terms.
These reductions lead respectively to unified product, double cross product,
semi-direct product, cocycle extension, and direct product $3$-Lie algebroids.

\noindent\textbf{Unified Product $3$-Lie Algebroids.}
If one of the twisted cocycle terms is zero, that is, if either $\psi$ or
$\zeta$ vanishes, then the bicocycle double cross product $3$-Lie algebroid
reduces to a unified product $3$-Lie algebroid.

Let us first assume that $\psi=0$. In this case, the operation
$[\bullet,\bullet,\bullet]_{\mathcal{B}}$ in \eqref{3-Lie-B} becomes closed,
thereby forming a genuine $3$-bracket on $\mathcal{B}$. Thus, the
$3$-anchored bundle $\mathcal{B}$ becomes a $3$-Lie subalgebroid of the product
space $\mathcal{A}\times\mathcal{B}$. Moreover, in this case the map $\rho$ in
\eqref{sennur3} becomes a genuine action. We denote the corresponding unified
product $3$-Lie algebroid by
\begin{equation}
\mathcal{A}_{\zeta}\bowtie\mathcal{B}.
\end{equation}
The $3$-Lie bracket \eqref{BDCP-3-Lie-bra}--\eqref{BDCP-3-Lie-bra+} reduces to
\begin{equation}
\begin{split}
&[(X_1,\eta_1),(X_2,\eta_2),(X_3,\eta_3)]_{{}_{\zeta}\bowtie}
=
(\widetilde{X},\widetilde{\eta}),
\\
&\widetilde{X}
=
\phi(X_1,X_2,X_3)
+
\rho(\eta_2,\eta_3,X_1)
-
\rho(\eta_1,\eta_3,X_2)
+
\rho(\eta_1,\eta_2,X_3),
\\
&\widetilde{\eta}
=
\theta(\eta_1,\eta_2,\eta_3)
+
\zeta(X_1,X_2,X_3)
+
\sigma(\eta_3,X_1,X_2)
-
\sigma(\eta_2,X_1,X_3)
+
\sigma(\eta_1,X_2,X_3).
\end{split}
\end{equation}

On the other hand, if $\zeta=0$, then the operation
$[\bullet,\bullet,\bullet]_{\mathcal{A}}$ in \eqref{3-Lie-A} becomes closed,
and hence it defines a genuine $3$-bracket on $\mathcal{A}$. Consequently,
$\mathcal{A}$ becomes a $3$-Lie subalgebroid of
$\mathcal{A}\times\mathcal{B}$. In this case, the map $\sigma$ in
\eqref{sennur3} becomes a genuine action. We denote this unified product
$3$-Lie algebroid by
\begin{equation}
\mathcal{A}\bowtie_{\psi}\mathcal{B}.
\end{equation}
The $3$-Lie bracket \eqref{BDCP-3-Lie-bra}--\eqref{BDCP-3-Lie-bra+} reduces to
\begin{equation}
\begin{split}
&[(X_1,\eta_1),(X_2,\eta_2),(X_3,\eta_3)]_{\bowtie_{\psi}}
=
(\widetilde{X},\widetilde{\eta}),
\\
&\widetilde{X}
=
\phi(X_1,X_2,X_3)
+
\psi(\eta_1,\eta_2,\eta_3)
+
\rho(\eta_2,\eta_3,X_1)
-
\rho(\eta_1,\eta_3,X_2)
+
\rho(\eta_1,\eta_2,X_3),
\\
&\widetilde{\eta}
=
\theta(\eta_1,\eta_2,\eta_3)
+
\sigma(\eta_3,X_1,X_2)
-
\sigma(\eta_2,X_1,X_3)
+
\sigma(\eta_1,X_2,X_3).
\end{split}
\end{equation}

\noindent\textbf{Double Cross and Semi-Direct Products.}
If both twisted cocycles are trivial, that is, if $\zeta=0$ and $\psi=0$, then
we arrive at the double cross product, or matched pair, $3$-Lie algebroid,
denoted by
\begin{equation}
\mathcal{A}\bowtie\mathcal{B}.
\end{equation}
In this case, the maps $\rho$ and $\sigma$ in \eqref{sennur3} are genuine
actions, and both $\mathcal{A}$ and $\mathcal{B}$ are $3$-Lie subalgebroids of
$\mathcal{A}\times\mathcal{B}$. The double cross product $3$-Lie bracket is
given by
\begin{equation}\label{DCP-brac}
\begin{split}
&[(X_1,\eta_1),(X_2,\eta_2),(X_3,\eta_3)]_{\bowtie}
=
(\widetilde{X},\widetilde{\eta}),
\\
&\widetilde{X}
=
\phi(X_1,X_2,X_3)
+
\rho(\eta_2,\eta_3,X_1)
-
\rho(\eta_1,\eta_3,X_2)
+
\rho(\eta_1,\eta_2,X_3),
\\
&\widetilde{\eta}
=
\theta(\eta_1,\eta_2,\eta_3)
+
\sigma(\eta_3,X_1,X_2)
-
\sigma(\eta_2,X_1,X_3)
+
\sigma(\eta_1,X_2,X_3).
\end{split}
\end{equation}

The semi-direct product $3$-Lie algebroids are obtained as one-sided special
cases of the double cross product. More precisely, they arise when one of the
actions $\rho$ or $\sigma$ in \eqref{DCP-brac} vanishes identically.

If $\rho=0$, then only the right action of $\mathcal{A}$ on $\mathcal{B}$
remains, while there is no action of $\mathcal{B}$ on $\mathcal{A}$. We denote
this semi-direct product by
\begin{equation}
\mathcal{A}\ltimes\mathcal{B}.
\end{equation}
The corresponding semi-direct product $3$-Lie bracket is
\begin{equation}\label{SD-brac}
\begin{split}
&[(X_1,\eta_1),(X_2,\eta_2),(X_3,\eta_3)]_{\ltimes}
=
(\widetilde{X},\widetilde{\eta}),
\\
&\widetilde{X}
=
\phi(X_1,X_2,X_3),
\\
&\widetilde{\eta}
=
\theta(\eta_1,\eta_2,\eta_3)
+
\sigma(\eta_3,X_1,X_2)
-
\sigma(\eta_2,X_1,X_3)
+
\sigma(\eta_1,X_2,X_3).
\end{split}
\end{equation}

If, on the other hand, $\sigma=0$, then only the left action of $\mathcal{B}$
on $\mathcal{A}$ remains, while there is no action of $\mathcal{A}$ on
$\mathcal{B}$. We denote this semi-direct product by
\begin{equation}
\mathcal{A}\rtimes\mathcal{B}.
\end{equation}
The corresponding semi-direct product $3$-Lie bracket is
\begin{equation}\label{SD-brac-}
\begin{split}
&[(X_1,\eta_1),(X_2,\eta_2),(X_3,\eta_3)]_{\rtimes}
=
(\widetilde{X},\widetilde{\eta}),
\\
&\widetilde{X}
=
\phi(X_1,X_2,X_3)
+
\rho(\eta_2,\eta_3,X_1)
-
\rho(\eta_1,\eta_3,X_2)
+
\rho(\eta_1,\eta_2,X_3),
\\
&\widetilde{\eta}
=
\theta(\eta_1,\eta_2,\eta_3).
\end{split}
\end{equation}
The discussion above may be summarized in the following statement.

\begin{lemma}\label{lem-special-cases}
The bicocycle double cross product $3$-Lie algebroid
$\mathcal{A}_{\zeta}\bowtie_{\psi}\mathcal{B}$ reduces to unified products,
double cross products, semi-direct products, cocycle extensions, and direct
products by imposing the corresponding vanishing conditions on the twisted
cocycles $\zeta,\psi$ and on the action terms $\rho,\sigma$.
\end{lemma}

The following diagram summarizes these reductions. In the diagram, DCP stands
for ``Double Cross Product'':
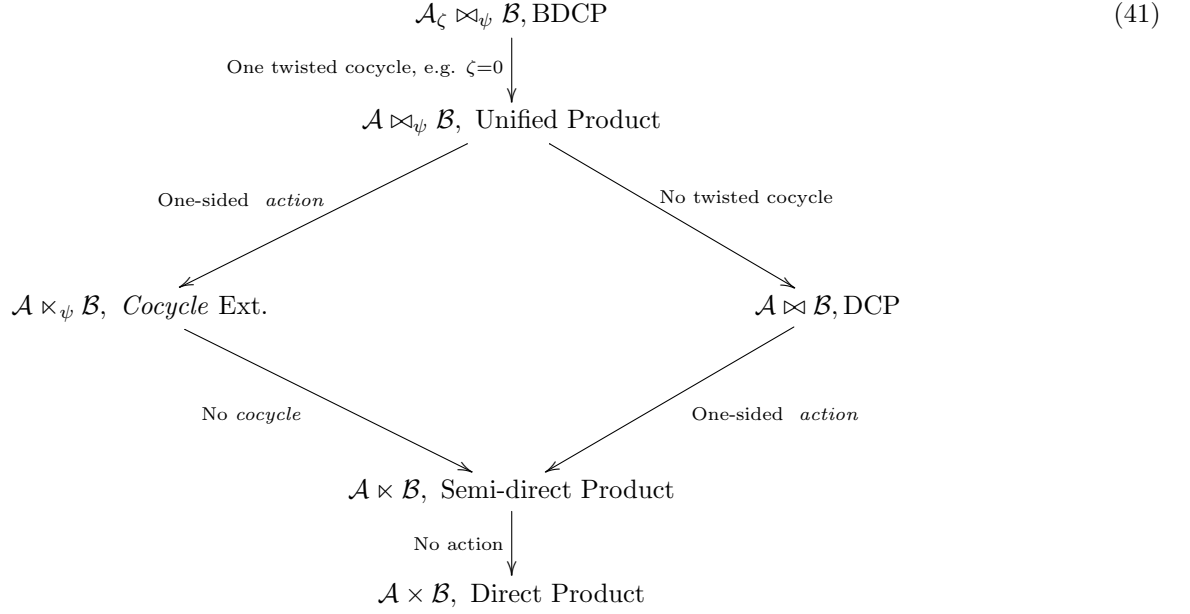
\begin{equation}\label{Diagram}
	\xymatrix{ 
 & \mathcal{A}_\zeta \bowtie _\psi\mathcal{B}, \text{BDCP}\ar[d]_{\text{One twisted cocycle, e.g. }\zeta=0}
 \\ &  \mathcal{A}\bowtie _\psi\mathcal{B}, \text{ Unified Product} \ar[ddl]_{\text{One-sided } \textit{  action}\qquad}\ar[ddr]^{\quad\text{No twisted cocycle}}\\
 & & & \\  \mathcal{A}\ltimes _\psi\mathcal{B},
		\text{ \textit{Cocycle} Ext.} 
		\ar[ddr]_{\text{No \textit{cocycle} } \qquad} &     &\mathcal{A}\bowtie\mathcal{B}, \text{DCP} \ar[ddl]^{\qquad\quad\text{One-sided } \textit{  action}\qquad}\\
  & & & \\
 & \mathcal{A}\ltimes \mathcal{B},
 \text{ Semi-direct Product} \ar[d]_{\text{No action}}
 \\
  & \mathcal{A}\times \mathcal{B}, \text{ Direct Product} 
 }
\end{equation}

\section{Decouplings in the BDCP Strategy}
Let us now examine the bicocycle double cross product strategy from the
opposite, namely decoupling, point of view. Consider a $3$-Lie algebroid
$\mathcal{M}$ and a subbundle $\mathcal{A}\subset\mathcal{M}$, not necessarily
a $3$-Lie subalgebroid. The $3$-Lie bracket on $\mathcal{M}$ induces
$3$-\textit{brackets} on both $\mathcal{A}$ and on a complementary bundle
$\mathcal{B}$. In general, these induced operations on $\mathcal{A}$ and
$\mathcal{B}$ need not be closed. The failure of closure gives rise to two
twisted cocycle terms: one describes how three elements of $\mathcal{A}$ yield
an element of $\mathcal{B}$, and the other describes how three elements of
$\mathcal{B}$ yield an element of $\mathcal{A}$. The mixed $3$-Lie brackets of
the form $[\mathcal{B},\mathcal{B},\mathcal{A}]$ and
$[\mathcal{B},\mathcal{A},\mathcal{A}]$ encode the mutual \textit{actions}.

In contrast to the complexity of the compatibility conditions
\eqref{doga1}, \eqref{pelin1}--\eqref{pelin5}, and \eqref{doga2}, the
bicocycle double cross product construction provides a flexible algebraic
language for decomposing a $3$-Lie algebroid. More precisely, assuming a mild
abuse of notation, we write
\begin{equation} 
X=(X,0)\in\Gamma(\mathcal{A}\times\mathcal{B}),
\qquad
\eta=(0,\eta)\in\Gamma(\mathcal{A}\times\mathcal{B}),
\end{equation}
and state the following result.

\begin{proposition}\label{mezgi+}  
Let
$(\mathcal{M},\upsilon,M,\mathfrak{a},[\bullet,\bullet,\bullet]_{\mathcal{M}})$
be a $3$-Lie algebroid, and assume that
$\mathcal{M}\cong\mathcal{A}\times\mathcal{B}$ as vector bundles, for two
complementary $3$-anchored bundles
\begin{equation}
(\mathcal{A},\tau=\upsilon\big|_{\mathcal{A}},M),
\qquad
(\mathcal{B},\kappa=\upsilon\big|_{\mathcal{B}},M).
\end{equation}
Assume further that the mixed brackets satisfy
\begin{equation}
[\eta_1,\eta_2,X]_{\mathcal{M}}\in\Gamma(\mathcal{A}),
\qquad
[\eta,X_1,X_2]_{\mathcal{M}}\in\Gamma(\mathcal{B}).
\end{equation}
Then $\mathcal{M}$ can be realized as the bicocycle double cross product of
$\mathcal{A}$ and $\mathcal{B}$, where the maps
\eqref{sennur2}, \eqref{sennur1}, and \eqref{sennur3} are obtained from the
decomposition
\begin{equation}\label{ezg}
\begin{split}
[X_1,X_2,X_3]_{\mathcal{M}}
&=
\phi(X_1,X_2,X_3)
+
\zeta(X_1,X_2,X_3),
\\
[\eta_1,\eta_2,\eta_3]_{\mathcal{M}}
&=
\psi(\eta_1,\eta_2,\eta_3)
+
\theta(\eta_1,\eta_2,\eta_3),
\\
[\eta_1,\eta_2,X]_{\mathcal{M}}
&=
\rho(\eta_1,\eta_2,X),
\\
[\eta,X_1,X_2]_{\mathcal{M}}
&=
\sigma(\eta,X_1,X_2).
\end{split}
\end{equation}
\end{proposition}

By considering the Whitney sum
$\mathcal{M}=\mathcal{A}\times\mathcal{B}$, we may list the possible
decomposition patterns from the simplest to the most general. We begin with the
case in which $\mathcal{A}$ is a $3$-Lie subalgebroid.

\begin{itemize}
    \item In the direct product case $\mathcal{A}\times\mathcal{B}$, the
    induced $3$-brackets on both $\mathcal{A}$ and $\mathcal{B}$ are closed.
    Therefore, no twisted cocycle terms appear, and the mixed $3$-brackets are
    identically zero. This corresponds to the case of trivial actions.

    \item In the semi-direct product case, there are still no twisted cocycle
    terms, since the induced $3$-brackets on the individual components are
    closed. However, the mixed brackets allow a one-sided action. In terms of
    \eqref{ezg}, this means that one of the maps $\rho$ or $\sigma$ is zero.
    These cases are denoted by $\mathcal{A}\ltimes\mathcal{B}$ and
    $\mathcal{A}\rtimes\mathcal{B}$, respectively.

    \item In a cocycle extension, a twisted cocycle term appears in one
    component while only a one-sided action is present. For example, the term
    $\psi$ records the component in $\mathcal{A}$ produced by the $3$-bracket
    of three elements of $\mathcal{B}$. We denote such a cocycle extension by
    $\mathcal{A}\ltimes_{\psi}\mathcal{B}$.

    \item In the double cross product, or matched pair, construction, both
    induced $3$-brackets are closed, so no twisted cocycle terms occur. However,
    the mixed brackets may have nontrivial components corresponding to mutual
    actions. We denote this case by $\mathcal{A}\bowtie\mathcal{B}$.

    \item In the unified product case, one twisted cocycle term is allowed
    together with mutual action terms. For instance, when the $3$-bracket of
    three elements of $\mathcal{B}$ has a nonzero component in $\mathcal{A}$,
    the cocycle $\psi$ appears. We denote this case by
    $\mathcal{A}\bowtie_{\psi}\mathcal{B}$.

    \item The most general case is obtained by allowing both twisted cocycle
    terms $\zeta$ and $\psi$, together with the two mutual action terms $\rho$
    and $\sigma$. This gives the bicocycle double cross product $3$-Lie
    algebroid, denoted by
    $\mathcal{A}_{\zeta}\bowtie_{\psi}\mathcal{B}$.
\end{itemize}

\section{BDCP 3-Lie Algebras}

As an immediate algebraic counterpart of the construction above, let
$\mathcal{A}$ and $\mathcal{B}$ be two vector spaces. We consider the
skew-symmetric trilinear maps
\begin{equation}\label{BDCP-3Lie-algebra-data-1}
\phi:\wedge^3\mathcal{A}\longrightarrow\mathcal{A},
\qquad
\zeta:\wedge^3\mathcal{A}\longrightarrow\mathcal{B},
\end{equation}
and
\begin{equation}\label{BDCP-3Lie-algebra-data-2}
\psi:\wedge^3\mathcal{B}\longrightarrow\mathcal{A},
\qquad
\theta:\wedge^3\mathcal{B}\longrightarrow\mathcal{B}.
\end{equation}
We also consider the action-type maps
\begin{equation}\label{BDCP-3Lie-algebra-data-3}
\rho:\wedge^2\mathcal{B}\otimes\mathcal{A}\longrightarrow\mathcal{A},
\qquad
\sigma:\mathcal{B}\otimes\wedge^2\mathcal{A}\longrightarrow\mathcal{B}.
\end{equation}
Using these data, we define a skew-symmetric trilinear bracket on
$\mathcal{A}\times\mathcal{B}$ by
\begin{equation}\label{BDCP-3Lie-algebra-bracket}
[(X_1,\eta_1),(X_2,\eta_2),(X_3,\eta_3)]_{{}_{\zeta}\bowtie_{\psi}}
=
(\widetilde{X},\widetilde{\eta}),
\end{equation}
where
\begin{equation}\label{BDCP-3Lie-algebra-bracket-components}
\begin{split}
\widetilde{X}
&=
\phi(X_1,X_2,X_3)
+
\psi(\eta_1,\eta_2,\eta_3)
+
\rho(\eta_2,\eta_3,X_1)
-
\rho(\eta_1,\eta_3,X_2)
+
\rho(\eta_1,\eta_2,X_3),
\\
\widetilde{\eta}
&=
\theta(\eta_1,\eta_2,\eta_3)
+
\zeta(X_1,X_2,X_3)
+
\sigma(\eta_3,X_1,X_2)
-
\sigma(\eta_2,X_1,X_3)
+
\sigma(\eta_1,X_2,X_3).
\end{split}
\end{equation}

\begin{proposition}\label{prop-BDCP-3Lie-algebra}
Let $\mathcal{A}$ and $\mathcal{B}$ be two vector spaces equipped with the maps
introduced in
\eqref{BDCP-3Lie-algebra-data-1},
\eqref{BDCP-3Lie-algebra-data-2}, and
\eqref{BDCP-3Lie-algebra-data-3}. Then the bracket
\eqref{BDCP-3Lie-algebra-bracket}--\eqref{BDCP-3Lie-algebra-bracket-components}
defines a $3$-Lie algebra structure on $\mathcal{A}\times\mathcal{B}$ if and
only if the compatibility conditions
\eqref{pelin1}--\eqref{pelin5}
are satisfied in the purely algebraic sense.
\end{proposition}

Whenever these conditions are satisfied, the resulting $3$-Lie algebra
\begin{equation}
\mathcal{A}_{\zeta}\bowtie_{\psi}\mathcal{B}
\end{equation}
is called the bicocycle double cross product $3$-Lie algebra associated with
the data
\begin{equation}
(\phi,\zeta,\psi,\theta,\rho,\sigma).
\end{equation}

The algebraic construction above immediately produces $3$-Lie algebroids.
Indeed, if $\mathcal{V}$ is a $3$-Lie algebra, then the trivial vector bundle
$M\times\mathcal{V}\longrightarrow M$, equipped with the zero anchor and with
the fiberwise $3$-Lie bracket extended $C^\infty(M)$-multilinearly, is a
$3$-Lie algebroid. In this case the anchor compatibility and Leibniz-type
conditions are automatically satisfied. Thus, every BDCP $3$-Lie algebra gives
rise to a zero-anchor BDCP $3$-Lie algebroid over an arbitrary smooth manifold
$M$.

\begin{example}[A double cross product $3$-Lie algebra and its trivial
algebroid]
Let
\begin{equation}
\mathcal{A}
=
\operatorname{span}\{\mathbf{e}_1,\mathbf{e}_2\},
\qquad
\mathcal{B}
=
\operatorname{span}\{\mathbf{e}_3,\mathbf{e}_4\}.
\end{equation}
Since both vector spaces have dimension two, the intrinsic skew-symmetric
$3$-brackets on them vanish. Thus, we set
\begin{equation}
\phi=0,
\qquad
\theta=0,
\qquad
\zeta=0,
\qquad
\psi=0.
\end{equation}
Define the action-type maps by
\begin{equation}
\rho(\mathbf{e}_3,\mathbf{e}_4,\mathbf{e}_1)=\mathbf{e}_2,
\qquad
\rho(\mathbf{e}_3,\mathbf{e}_4,\mathbf{e}_2)=-\mathbf{e}_1,
\end{equation}
and
\begin{equation}
\sigma(\mathbf{e}_3,\mathbf{e}_1,\mathbf{e}_2)=\mathbf{e}_4,
\qquad
\sigma(\mathbf{e}_4,\mathbf{e}_1,\mathbf{e}_2)=-\mathbf{e}_3,
\end{equation}
and extend these maps by skew-symmetry.

The corresponding double cross product bracket on
$\mathcal{A}\times\mathcal{B}$ is determined by
\begin{equation}
[(\mathbf{e}_1,0),(\mathbf{e}_2,0),(0,\mathbf{e}_3)]_{\bowtie}
=
(0,\mathbf{e}_4),
\qquad
[(\mathbf{e}_1,0),(\mathbf{e}_2,0),(0,\mathbf{e}_4)]_{\bowtie}
=
-(0,\mathbf{e}_3),
\end{equation}
and
\begin{equation}
[(\mathbf{e}_1,0),(0,\mathbf{e}_3),(0,\mathbf{e}_4)]_{\bowtie}
=
(\mathbf{e}_2,0),
\qquad
[(\mathbf{e}_2,0),(0,\mathbf{e}_3),(0,\mathbf{e}_4)]_{\bowtie}
=
-(\mathbf{e}_1,0),
\end{equation}
together with skew-symmetric permutations. Equivalently, if
\begin{equation}
\mathbf{E}_1=(\mathbf{e}_1,0),
\quad
\mathbf{E}_2=(\mathbf{e}_2,0),
\quad
\mathbf{E}_3=(0,\mathbf{e}_3),
\quad
\mathbf{E}_4=(0,\mathbf{e}_4),
\end{equation}
then the nonzero brackets are
\begin{equation}
[\mathbf{E}_1,\mathbf{E}_2,\mathbf{E}_3]_{\bowtie}
=
\mathbf{E}_4,
\qquad
[\mathbf{E}_1,\mathbf{E}_2,\mathbf{E}_4]_{\bowtie}
=
-\mathbf{E}_3,
\end{equation}
and
\begin{equation}
[\mathbf{E}_1,\mathbf{E}_3,\mathbf{E}_4]_{\bowtie}
=
\mathbf{E}_2,
\qquad
[\mathbf{E}_2,\mathbf{E}_3,\mathbf{E}_4]_{\bowtie}
=
-\mathbf{E}_1.
\end{equation}
This is the standard four-dimensional simple $3$-Lie algebra. Therefore,
$\mathcal{A}\bowtie\mathcal{B}$ is a double cross product $3$-Lie algebra with
nontrivial mutual actions $\rho$ and $\sigma$.

Consequently, for every smooth manifold $M$, the trivial vector bundle
\begin{equation}
M\times(\mathcal{A}\times\mathcal{B})
\longrightarrow M,
\end{equation}
equipped with the zero anchor and the fiberwise bracket above, is a double cross
product $3$-Lie algebroid.
\end{example}

\begin{example}[A BDCP $3$-Lie algebra with all structure maps nontrivial]
Let
\begin{equation}
\mathcal{A}
=
\operatorname{span}
\{\mathbf{a}_1,\mathbf{a}_2,\mathbf{a}_3,\mathbf{p},\mathbf{q},\mathbf{r}\},
\qquad
\mathcal{B}
=
\operatorname{span}
\{\mathbf{b}_1,\mathbf{b}_2,\mathbf{b}_3,\mathbf{u},\mathbf{v},\mathbf{w}\}.
\end{equation}
Define
\begin{equation}
\phi(\mathbf{a}_1,\mathbf{a}_2,\mathbf{a}_3)=\mathbf{p},
\qquad
\zeta(\mathbf{a}_1,\mathbf{a}_2,\mathbf{a}_3)=\mathbf{u},
\end{equation}
and
\begin{equation}
\psi(\mathbf{b}_1,\mathbf{b}_2,\mathbf{b}_3)=\mathbf{q},
\qquad
\theta(\mathbf{b}_1,\mathbf{b}_2,\mathbf{b}_3)=\mathbf{v}.
\end{equation}
Moreover, define
\begin{equation}
\rho(\mathbf{b}_1,\mathbf{b}_2,\mathbf{a}_1)=\mathbf{r},
\qquad
\sigma(\mathbf{b}_1,\mathbf{a}_1,\mathbf{a}_2)=\mathbf{w}.
\end{equation}
All these maps are extended by skew-symmetry, and all other values are defined
to be zero.

The induced BDCP bracket on $\mathcal{A}\times\mathcal{B}$ has the following
nonzero values:
\begin{equation}
\begin{split}
&[(\mathbf{a}_1,0),(\mathbf{a}_2,0),(\mathbf{a}_3,0)]_{{}_{\zeta}\bowtie_{\psi}}
=
(\mathbf{p},\mathbf{u}),
\\
&[(0,\mathbf{b}_1),(0,\mathbf{b}_2),(0,\mathbf{b}_3)]_{{}_{\zeta}\bowtie_{\psi}}
=
(\mathbf{q},\mathbf{v}),
\\
&[(0,\mathbf{b}_1),(0,\mathbf{b}_2),(\mathbf{a}_1,0)]_{{}_{\zeta}\bowtie_{\psi}}
=
(\mathbf{r},0),
\\
&[(0,\mathbf{b}_1),(\mathbf{a}_1,0),(\mathbf{a}_2,0)]_{{}_{\zeta}\bowtie_{\psi}}
=
(0,\mathbf{w}),
\end{split}
\end{equation}
together with their skew-symmetric permutations.

The elements
\begin{equation}
(\mathbf{p},0),
\quad
(\mathbf{q},0),
\quad
(\mathbf{r},0),
\quad
(0,\mathbf{u}),
\quad
(0,\mathbf{v}),
\quad
(0,\mathbf{w})
\end{equation}
are central with respect to this bracket. Therefore, for all
$U_1,U_2,V_1,V_2,V_3$ in $\mathcal{A}\times\mathcal{B}$, we have
\begin{equation}
[U_1,U_2,[V_1,V_2,V_3]_{{}_{\zeta}\bowtie_{\psi}}]_{{}_{\zeta}\bowtie_{\psi}}
=
0.
\end{equation}
Each term on the right-hand side of the fundamental identity vanishes for the
same reason. Hence the fundamental identity is satisfied, and 
$\mathcal{A}_{\zeta}\bowtie_{\psi}\mathcal{B}$ 
is a BDCP $3$-Lie algebra. Moreover, all six structure maps  
$\phi, \zeta,  \psi, \theta, \rho$, and $ \sigma$ 
are nontrivial.

Consequently, for every smooth manifold $M$, the trivial vector bundle
\begin{equation}
M\times(\mathcal{A}\times\mathcal{B})
\longrightarrow M,
\end{equation}
equipped with the zero anchor and the fiberwise bracket above extended
$C^\infty(M)$-multilinearly, is a BDCP $3$-Lie algebroid with all structure maps
nontrivial except for the anchor.
\end{example}

 \section{Conclusion}

In this work, we have developed an algebraic framework for merging two
$3$-anchored bundles equipped with mutual \textit{actions} and two twisted
cocycle terms into a $3$-Lie algebroid. The main result, Proposition
\ref{thm-muratcan-off+}, provides necessary and sufficient compatibility
conditions under which the Whitney sum of two $3$-anchored bundles carries
a bicocycle double cross product $3$-Lie algebroid structure.

We have also recorded the purely algebraic counterpart of this construction.
More precisely, Proposition \ref{prop-BDCP-3Lie-algebra} shows that, by
removing the anchor compatibility and Leibniz-type conditions from the
$3$-Lie algebroid setting, one obtains the bicocycle double cross product
construction for $3$-Lie algebras. In this purely algebraic case, the remaining
compatibility conditions are precisely those arising from the fundamental
identity.

Furthermore, Lemma \ref{lem-special-cases} shows that the proposed construction
contains several important structures as special cases. These include unified
products, matched pairs or double cross products, semi-direct products, cocycle
extensions, and direct products. These reductions are also summarized in
Diagram \ref{Diagram}.

Finally, we have examined the converse problem from the decoupling point of
view. Proposition \ref{mezgi+} shows how a $3$-Lie algebroid decomposed into two
complementary $3$-anchored bundles gives rise to the algebraic data of the
bicocycle double cross product construction, namely componentwise brackets,
mutual action terms, and twisted cocycles.

Thus, the bicocycle double cross product provides a unifying algebraic scheme
for both constructing and decomposing $3$-Lie algebroids and $3$-Lie algebras.

\section*{Acknowledgment}

The authors acknowledge the use of GPT for grammatical checking and language
editing at the final stage of the manuscript preparation.

\bibliographystyle{abbrv}
\bibliography{references}

\begin{thebibliography}{10}

\bibitem{AgorMili11}
A.~L. Agore and G.~Militaru.
\newblock Extending structures {II}: {T}he quantum version.
\newblock {\em J. Algebra}, 336:321--341, 2011.

\bibitem{Agore-Lie}
A.~L. Agore and G.~Militaru.
\newblock Extending structures for {L}ie algebras.
\newblock {\em Monatsh. Math.}, 174(2):169--193, 2014.

\bibitem{AgorMili14-II}
A.~L. Agore and G.~Militaru.
\newblock Extending structures {I}: the level of groups.
\newblock {\em Algebr. Represent. Theory}, 17(3):831--848, 2014.

\bibitem{AgoreMilitaru-book}
A.~L. Agore and G.~Militaru.
\newblock {\em Extending {S}tructures: {F}undamentals and {A}pplications}.
\newblock CRC Press, Taylor \& Francis Group, 2019.
\newblock Chapman \& Hall/CRC Monographs and Research Notes in Mathematics.

\bibitem{AjChEsGrKlPa}
A.~Ajji, J.~Chaouki, O.~Esen, M.~Grmela, V.~Klika, and M.~Pavelka.
\newblock On geometry of multiscale mass action law and its fluctuations.
\newblock {\em Phys. D}, 445:Paper No. 133642, 22, 2023.

\bibitem{atecsli2025product}
B.~Ate{\c{s}}li, O.~Esen, and S.~S{\"u}tl{\"u}.
\newblock On product {L}ie algebroids, and collective motion.
\newblock {\em Proceedings of the Institute of Mathematics and Mechanics},
  51(2):205–244, 2025.

\bibitem{uccgun2024dynamics}
F.~{\c{C}}a{\u{g}}atay~U{\c{c}}gun, O.~Esen, and S.~S{\"u}tl{\"u}.
\newblock Dynamics over cocycle double cross-products.
\newblock {\em International Journal of Geometric Methods in Modern Physics},
  page 2550036, 2025.

\bibitem{CeHoMaRa98}
H.~Cendra, D.~D. Holm, J.~E. Marsden, and T.~S. Ratiu.
\newblock Lagrangian reduction, the {E}uler-{P}oincar\'{e} equations, and
  semidirect products.
\newblock In {\em Geometry of differential equations}, volume 186 of {\em Amer.
  Math. Soc. Transl. Ser. 2}, pages 1--25. Amer. Math. Soc., Providence, RI,
  1998.

\bibitem{EsGrMiGu19}
O.~Esen, M.~Grmela, H.~G\"{u}mral, and M.~Pavelka.
\newblock Lifts of symmetric tensors: fluids, plasma, and {G}rad hierarchy.
\newblock {\em Entropy}, 21(9):Paper No. 907, 33, 2019.

\bibitem{esen2021bicocycle}
O.~Esen, P.~Guha, and S.~S\"{u}tl\"{u}.
\newblock Bicocycle double cross constructions.
\newblock {\em J. Algebra Appl.}, 22(12):Paper No. 2350254, 36, 2023.

\bibitem{EsenKudeSutlu21}
O.~Esen, M.~Kudeyt, and S.~S\"{u}tl\"{u}.
\newblock Second order {L}agrangian dynamics on double cross product groups.
\newblock {\em J. Geom. Phys.}, 159:103934, 18, 2021.

\bibitem{esen2021extensions}
O.~Esen, G.~\"{O}zcan, and S.~S\"{u}tl\"{u}.
\newblock On extensions, {L}ie-{P}oisson systems, and dissipation.
\newblock {\em J. Lie Theory}, 32(2):327--382, 2022.

\bibitem{Esvd17}
O.~Esen, M.~Pavelka, and M.~Grmela.
\newblock Hamiltonian coupling of electromagnetic field and matter.
\newblock {\em Int. J. Adv. Eng. Sci. Appl. Math.}, 9(1):3--20, 2017.

\bibitem{esen2016hamiltonian}
O.~Esen and S.~S{\"u}tl{\"u}.
\newblock Hamiltonian dynamics on matched pairs.
\newblock {\em International Journal of Geometric Methods in Modern Physics},
  13(10):1650128, 2016.

\bibitem{esen2017lagrangian}
O.~Esen and S.~S{\"u}tl{\"u}.
\newblock Lagrangian dynamics on matched pairs.
\newblock {\em Journal of Geometry and Physics}, 111:142--157, 2017.

\bibitem{esen2018matched}
O.~Esen and S.~S{\"u}tl{\"u}.
\newblock Discrete dynamical systems over double cross-product lie groupoids.
\newblock {\em International Journal of Geometric Methods in Modern Physics},
  18(04):2150057, 2021.

\bibitem{EsSu21}
O.~Esen and S.~S\"{u}tl\"{u}.
\newblock Matched pair analysis of the {V}lasov plasma.
\newblock {\em J. Geom. Mech.}, 13(2):209--246, 2021.

\bibitem{GrabowskaAlg}
K.~Grabowska, P.~Urba\'{n}ski, and J.~Grabowski.
\newblock Geometrical mechanics on algebroids.
\newblock {\em Int. J. Geom. Methods Mod. Phys.}, 3(3):559--575, 2006.

\bibitem{GrabowskiMarmo2000}
J.~Grabowski and G.~Marmo.
\newblock On {F}ilippov algebroids and multiplicative {N}ambu-{P}oisson
  structures.
\newblock {\em Differential Geom. Appl.}, 12(1):35--50, 2000.

\bibitem{kuper83}
B.~A. Kupershmidt and T.~Ratiu.
\newblock Canonical maps between semidirect products with applications to
  elasticity and superfluids.
\newblock {\em Comm. Math. Phys.}, 90(2):235--250, 1983.

\bibitem{MackenzieDG}
K.~C.~H. Mackenzie.
\newblock {\em Lie groupoids and {L}ie algebroids in differential geometry},
  volume 124 of {\em London Mathematical Society Lecture Note Series}.
\newblock Cambridge University Press, Cambridge, 1987.

\bibitem{Mackenzie-book}
K.~C.~H. Mackenzie.
\newblock {\em General theory of {L}ie groupoids and {L}ie algebroids}, volume
  213 of {\em London Mathematical Society Lecture Note Series}.
\newblock Cambridge University Press, Cambridge, 2005.

\bibitem{Ma90}
S.~Majid.
\newblock Matched pairs of {L}ie groups associated to solutions of the
  {Y}ang-{B}axter equations.
\newblock {\em Pacific Journal of Mathematics}, 141(2):311--332, 1990.

\bibitem{majid1990physicsa}
S.~Majid.
\newblock Physics for algebraists: Non-commutative and non-cocommutative hopf
  algebras by a bicrossproduct construction.
\newblock {\em Journal of Algebra}, 130(1):17--64, 1990.

\bibitem{marle2014lie}
C.~Marle.
\newblock Lie, symplectic and {P}oisson groupoids and their {L}ie algebroids.
\newblock {\em arXiv preprint arXiv:1402.0059}, 2014.

\bibitem{marsden83b}
J.~E. Marsden, T.~Ratiu, and A.~Weinstein.
\newblock Reduction and {H}amiltonian structures on duals of semidirect product
  {L}ie algebras.
\newblock In {\em Fluids and plasmas: geometry and dynamics ({B}oulder,
  {C}olo., 1983)}, volume~28 of {\em Contemp. Math.}, pages 55--100. Amer.
  Math. Soc., Providence, RI, 1984.

\bibitem{marsden1984semidirect}
J.~E. Marsden, T.~Ra{\c{t}}iu, and A.~Weinstein.
\newblock Semidirect products and reduction in mechanics.
\newblock {\em Transactions of the american mathematical society},
  281(1):147--177, 1984.

\bibitem{martinez2001}
E.~Mart\'{\i}nez.
\newblock Lagrangian mechanics on {L}ie algebroids.
\newblock {\em Acta Appl. Math.}, 67(3):295--320, 2001.

\bibitem{Mishra}
S.~K. Mishra, G.~Mukherjee, and A.~Naolekar.
\newblock Cohomology and deformations of {F}ilippov algebroids.
\newblock {\em Proc. Indian Acad. Sci. Math. Sci.}, 132(1):Paper No. 2, 25,
  2022.

\bibitem{Mokr97}
T.~Mokri.
\newblock Matched pairs of {L}ie algebroids.
\newblock {\em Glasgow Math. J.}, 39(2):167--181, 1997.

\bibitem{Para67}
J.~Pradines.
\newblock Th\'{e}orie de {L}ie pour les groupo\"{\i}des diff\'{e}rentiables.
  {C}alcul diff\'{e}renetiel dans la cat\'{e}gorie des groupo\"{\i}des
  infinit\'{e}simaux.
\newblock {\em C. R. Acad. Sci. Paris S\'{e}r. A-B}, 264:A245--A248, 1967.

\bibitem{Vallejo}
J.~A. Vallejo.
\newblock Nambu-{P}oisson manifolds and associated {$n$}-ary {L}ie algebroids.
\newblock {\em J. Phys. A}, 34(13):2867--2881, 2001.

\bibitem{WeinLag}
A.~Weinstein.
\newblock {L}agrangian mechanics and groupoids.
\newblock In {\em Mechanics day ({W}aterloo, {ON}, 1992)}, volume~7 of {\em
  Fields Inst. Commun.}, pages 207--231. Amer. Math. Soc., Providence, RI,
  1996.

\end{thebibliography}

\end{document}